\documentclass{amsart}

\def\R{{\mathbb {R}}}
\def\N{{\mathbb {N}}}

\newtheorem{teo}{Theorem}
\newtheorem{lema}{Lemma}

\newtheorem{prop}{Proposition}

\theoremstyle{remark}
\newtheorem{remark}{Remark}

\begin{document}

\title[Multiple solutions with critical growth]{Multiple
solutions for the $p(x)-$laplace operator with critical growth}

\author[ A. Silva]
{ Anal\'{\i}a Silva}

\address{Anal\'ia Silva\hfill\break\indent
Departamento  de Matem\'atica, FCEyN \hfill\break\indent UBA (1428) Buenos Aires, Argentina.}
\email{{\tt asilva@dm.uba.ar}}

\thanks{Supported by Universidad de Buenos Aires under grant X078 and X837,
by ANPCyT PICT No. 2006-290 and CONICET (Argentina) PIP 5477 and PIP 5478/1438.Anal\'ia Silva is a fellow of CONICET.}

\begin{abstract}
The aim of this paper is to extend previous results regarding the multiplicity of solutions for quasilinear elliptic problems with critical growth to the variable exponent case.

We prove, in the spirit of \cite{DPFBS}, the existence of at least three nontrivial solutions to the following quasilinear elliptic equation $-\Delta_{p(x)} u = |u|^{q(x)-2}u +\lambda f(x,u)$ in a smooth bounded domain $\Omega$ of $\R^N$ with homogeneous Dirichlet boundary conditions on $\partial\Omega$. We assume that $\{q(x)=p^*(x)\}\not=\emptyset$, where $p^*(x)=Np(x)/(N-p(x))$ is the critical Sobolev exponent for variable exponents and $\Delta_{p(x)} u = \mbox{div}(|\nabla u|^{p(x)-2}\nabla u)$ is the $p(x)-$laplacian. The proof is based on variational arguments and the extension of  concentration compactness method for variable exponent spaces.
\end{abstract}

\maketitle

\section{Introduction.}

Let us consider the following nonlinear elliptic problem:
\begin{equation}\tag{P}\label{1.1}
\begin{cases}
-\Delta_{p(x)} u = |u|^{q(x)-2}u + \lambda f(x,u) & \mbox{in } \Omega\\
u = 0 & \mbox{on } \partial\Omega,
\end{cases}
\end{equation}
where $\Omega$ is a bounded smooth domain in $\R^N$, $\Delta_{p(x)} u =
\mbox{div}(|\nabla u|^{p(x)-2}\nabla u)$ is the $p(x)-$laplacian, $1<p(x)<N$. On the exponent $q(x)$ we assume that is critical in the sense that $\{q(x)=p^*(x)\}\not=\emptyset$,
where $p^*(x)=Np(x)/(N-p(x))$ is the critical exponent in the Sobolev embedding, $\lambda$ is a positive parameter and the nonlinear term $f$ is a subcritical perturbation with some precise assumptions that we state below.

The purpose of this paper, is to extend the results obtained in \cite{DPFBS} where the same problem but with constant $p$ was treated. Namely, in \cite{DPFBS}, problem \eqref{1.1} was analyzed in the case $p(x)\equiv p$ constant and $q(x)\equiv p^*$.

To be more precise, the result in \cite{DPFBS} prove of the existence of at least three nontrivial solutions for \eqref{1.1}, one positive, one negative and one that changes sign, under adequate assumptions on the source term $f$ and the parameter $\lambda$.

The method in the proof used in \cite{DPFBS} consists on restricting the functional associated to \eqref{1.1} to three different Banach manifolds, one consisting on positive functions, one consisting on negative functions  and the third one consisting on sign-changing functions, all of them under a normalization condition, Then, by means of a suitable version of the Mountain Pass Theorem due to Schwartz \cite{S} and the concentration-compactness principle of P.L. Lions \cite{Lions} the authors can prove the existence of a critical point of each restricted functional and, finally, the authors were able to prove that critical points of each restricted functional are critical points of the unrestricted one.

This method was introduced by M. Struwe \cite{St} where the subcritical case (in the sense of the Sobolev embeddigs) for the $p-$Laplacian was treated. A related result for the $p-$Laplacian under nonlinear boundary condition can be found in \cite{FB}.

Also, a similar problem in the case of the $p(x)-$Laplacian, but with subcritical nonlinearities was analyzed in \cite{DU}.

In all the above mentioned works, the main feature on the nonlinear term $f$ is that no oddness condition is imposed

Very little is known about critical growth nonlinearities for variable exponent problems, since one of the main techniques used in order to deal with such issues is the concentration-compactness principle. This result was recently obtained for the variable exponent case independently in \cite{FBS} and \cite{Fu}.
In both of these papers the proof are similar and both relates to that of the original proof of P.L. Lions. However, the arguments in \cite{FBS} are a little more subtle and allow the authors to deal with the case where the exponent $q(x)$ is critical only in some part of the domain, while the results in \cite{Fu} requires $q(x)$ to be identically $p^*(x)$. So we will rely on the concentration-compactness principle proved in \cite{FBS} in this work.

The use of the concentration compactness method to deal with the $p-$Laplacian has been used by so many authors before that is almost impossible to give a complete list of contributions. However we want to refer to the work of J. Garc\'{\i}a Azorero and I. Peral in \cite{GAP} from where we borrow some ideas.

Throughout this work, by (weak) solutions of \eqref{1.1} we understand critical points of the associated energy functional acting on the Sobolev space $W_0^{1,p(x)}(\Omega)$:
\begin{equation}\label{Phi}
\Phi(v)=\int_{\Omega}\frac{1}{p(x)} |\nabla v|^{p(x)} \, dx -\int_{\Omega}\frac{1}{q(x)}
|u|^{q(x)}\, dx - \lambda \int_{\Omega} F(x,v)\, dx,
\end{equation}
where $F(x,u) = \int_0^u f(x,z)\, dz$.

To end this introduction, let us comment on different applications where the $p(x)-$Laplacian has appeared.

Up to our knowledge there are two main fields where the $p(x)-$Laplacian have been proved to be extremely usefull in applications:
\begin{itemize}
\item Image Processing
\item Electrorheological Fluids
\end{itemize}

For instance, Y. Chen, S. Levin and R. Rao \cite{CLR} proposed the following model in image processing
$$
E(u)=\int_{\Omega}\frac{|\nabla u(x)|^{p(x)}}{p(x)}+ f(|u(x)-I(x)|)\, dx \to \mbox{min}
$$
where $p(x)$ is a function varying between $1$ and $2$ and $f$ is a convex function.

In their application, they chose $p(x)$ close to 1 where there is likely to be edges and close to 2 where it is likely not to be edges.

The electrorheological fluids application is much more developed and we refer to the monograph by M. Ru{\v{z}}i{\v{c}}ka, \cite{Ru}, and its references.


\section{Assumptions and statement of the results.}

Throughout this paper the following notation will be used: Given $q\colon \Omega\to\R$ bounded, we denote
$$
q^+ := \sup_\Omega q(x), \qquad q^- := \inf_\Omega q(x).
$$
The precise assumptions on the source term $f$ are as follows:

\begin{enumerate}
\item[(F1)] $f:\Omega\times\R\to\R$, is a measurable function with respect to the first argument and continuously differentiable with respect to the second argument for almost every $x\in\Omega$. Moreover, $f(x,0)=0$ for every $x\in\Omega$.


\item[(F2)] There exist constants $c_1> 1/(q^--1)$, $c_2\in(p^+,q^-)$, $0<c_3<c_4$, such that for any $u\in L^q(\Omega)$ and $p^-\leq p^+<r^-\leq r^+<q^-\leq q^+$.
\begin{align*}
c_3\rho_r(u)&\leq c_2\int_\Omega F(x,u)\,dx\leq\int_\Omega f(x,u)u\,dx\\
&\leq c_1\int_\Omega f_u(x,u)u^2\,dx\leq c_4\rho_r(u)
\end{align*}
Where $\rho_r(u):=\int_\Omega|u|^{r(x)}\,dx$
\end{enumerate}

\begin{remark}
Observe that this set of hypotheses on the nonlinear term $f$ are weaker than the ones considered by \cite{S}.
\end{remark}

\begin{remark}
We exhibit now one example of nonlinearities that fulfill all of our hypotheses.
 $f(x,u)=|u|^{r(x)-2}u+|u_+|^{s(x)-2}u_+$, if $s(x)<r(x)$ , $q^--1>s^->p^+$.

Hypotheses (F1)--(F2) are clearly satisfied.

\end{remark}

So the main result of the paper reads:

\begin{teo} Under assumptions {\em (F1)--(F2)}, there
exist $\lambda^*>0$ depending only on $n, p, q$ and the constant $c_3$ in {\em (F2)}, such that for every $\lambda>\lambda^*$, there exists three different, nontrivial, (weak) solutions of problem \eqref{1.1}. Moreover these solutions are, one positive, one negative and the other one has non-constant sign.
\end{teo}

\section{Results on variable exponent Sobolev spaces}

The variable exponent Lebesgue space $L^{p(x)}(\Omega)$ is defined by
$$
L^{p(x)}(\Omega) = \left\{u\in L^1_{\text{loc}}(\Omega) \colon \int_\Omega|u(x)|^{p(x)}\,dx<\infty\right\}.
$$
This space is endowed with the norm
$$
\|u\|_{L^{p(x)}(\Omega)}=\inf\left\{\lambda>0:\int_\Omega\left|\frac{u(x)}{\lambda}\right|^{p(x)}\,dx\leq 1\right\}
$$
The variable exponent Sobolev space $W^{1,p(x)}(\Omega)$ is defined by
$$
W^{1,p(x)}(\Omega) = \{u\in W^{1,1}_{\text{loc}}(\Omega) \colon u\in L^{p(x)}(\Omega) \mbox{ and } |\nabla u |\in  L^{p(x)}(\Omega)\}.
$$
The corresponding norm for this space is
$$
\|u\|_{W^{1,p(x)}(\Omega)}=\|u\|_{L^{p(x)}(\Omega)}+\| |\nabla u| \|_{L^{p(x)}(\Omega)}
$$
Define $W^{1,p(x)}_0(\Omega)$ as the closure of $C_0^\infty(\Omega)$ with respect to the $W^{1,p(x)}(\Omega)$ norm. The spaces $L^{p(x)}(\Omega)$, $W^{1,p(x)}(\Omega)$ and $W^{1,p(x)}_0(\Omega)$ are separable and reflexive Banach spaces when $1<\inf_\Omega p \le \sup_\Omega p <\infty$.

As usual, we denote $p'(x) = p(x)/(p(x)-1)$ the conjugate exponent of $p(x)$.

Define
$$
p^*(x)=\begin{cases}
\frac{Np(x)}{N-p(x)} & \mbox{ if } p(x)<N\\
\infty & \mbox{ if } p(x)\geq N
\end{cases}
$$

The following results are proved in \cite{Fan}
\begin{prop}[H\"older-type inequality]\label{Holder}
Let $f\in L^{p(x)}(\Omega)$ and $g\in L^{p'(x)}(\Omega)$. Then the following inequality holds
$$
\int_\Omega|f(x)g(x)|\,dx\leq C_p\|f\|_{L^{p(x)}(\Omega)}\|g\|_{L^{p'(x)}(\Omega)}
$$
\end{prop}

\begin{prop}[Sobolev embedding]\label{embedding}
Let $p, q\in C(\overline{\Omega})$ be such that $1\leq q(x)\le p^*(x)$ for all $x\in\overline{\Omega}$. Assume moreover that the functions $p$ and $q$ are log-H\"older continuous. Then there is a continuous embedding
$$
W^{1,p(x)}(\Omega)\hookrightarrow L^{q(x)}(\Omega).
$$
Moreover, if $\inf_{\Omega} (p^*-q)>0$ then, the embedding is compact.
\end{prop}

\begin{prop}[Poincar\'e inequality]\label{Poincare}
There is a constant $C>0$, such that
$$
\|u\|_{L^{p(x)}(\Omega)}\leq C\| |\nabla u| \|_{L^{p(x)}(\Omega)},
$$
for all $u\in W^{1,p(x)}_0(\Omega)$.
\end{prop}

\begin{remark}
By Proposition \ref{Poincare}, we know that $\| |\nabla u| \|_{L^{p(x)}(\Omega)}$ and $\|u\|_{W^{1,p(x)}(\Omega)}$ are equivalent norms on $W_0^{1,p(x)}(\Omega)$.
\end{remark}

\section{Proof of Theorem 1.}

The proof uses the same approach as in \cite{St}. That is, we will construct three disjoint sets $K_i\ne\emptyset$ not containing $0$ such that $\Phi$ has a critical point in $K_i$. These sets will be subsets of $C^{1}-$manifolds $M_i\subset W^{1,p(x)}(\Omega)$ that will be constructed by imposing a sign restriction and a normalizing condition.

In fact, let
\begin{align*}
& J(v) =    \int_{\Omega}|\nabla v|^{p(x)}-|v|^{q(x)}dx\\
& M_1 = \left\{u \in W^{1,p(x)}_0(\Omega)\colon \int_{\Omega} u_+>0 \mbox{ and }
J(u_+) = \int_{\Omega}\lambda f(x,u) u_+
dx\right\},\\
& M_2 = \left\{u \in W^{1,p(x)}_0(\Omega)\colon \int_\Omega u_->0 \mbox{ and }
J(u_-) = -\int_\Omega\lambda f(x,u) u_- dx\right\},\\
& M_3 = M_1\cap M_2.
\end{align*}
where $u_+ = \max\{u,0\}$, $u_-=\max\{-u,0\}$ are the positive and negative parts of $u$.

Finally we define
\begin{align*}
& K_1 = \{ u\in M_1\ |\ u\ge 0\},\\
& K_2 = \{ u\in M_2\ |\ u\le 0\},\\
& K_3 = M_3.
\end{align*}

First, we need a Lemma to show that these sets are nonempty and, moreover, give some properties that will be useful in the proof of the result.

\begin{lema}\label{tlambda}
For every $w_0\in W^{1,p(x)}_0(\Omega)$, $w_0>0$ ($w_0<0$), there exists $t_\lambda>0$ such that $t_\lambda w_0\in M_1 (\in M_2)$. Moreover, $\lim_{\lambda\to\infty} t_\lambda = 0$.

As a consequence, given $w_0, w_1\in W^{1,p(x)}_0(\Omega)$, $w_0>0$, $w_1<0$, with disjoint supports, there exists $\bar t_\lambda, \underbar t_\lambda>0$ such that $\bar t_\lambda w_0 + \underbar t_\lambda w_1 \in M_3$. Moreover $\bar t_\lambda, \underbar t_\lambda \to 0$ as $\lambda\to\infty$.
\end{lema}

\begin{proof}
We prove the lemma for $M_1$, the other cases being similar.

For $w\in W^{1,p(x)}_0(\Omega)$, $w \geq 0$, we consider the functional
$$
\varphi_1(w)= \int_{\Omega} |\nabla w|^{p(x)} -|w|^{q(x)} - \lambda f(x,w)w \, dx.
$$
Given $w_0>0$, in order to prove the lemma, we must show that $\varphi_1(t_\lambda w_0)=0$ for some $t_\lambda>0$. Using hypothesis (F2), if $t<1$, we have that:
$$
\varphi_1(tw_0) \geq A t^{p^+} - B t^{q^-} - \lambda c_4 C t^{r^-}
$$
and
$$
\varphi_1(tw_0) \leq A t^{p^-} - B t^{q^+} - \lambda c_3 C t^{r^+},
$$
where the coefficients $A$, $B$ and $C$ are given by:
$$
A = \int_{\Omega} |\nabla w_0|^{p(x)} \,dx, \quad B= \int_{\Omega} |w_0|^{q(x)} \,dx,
\quad C = \int_{\Omega} |w_0|^{r(x)} \, dx.
$$
Since $p^-\leq p^+<r^-\leq r^+<q^-\leq q^+$ it follows that $\varphi_1(tw_0)$ is positive for $t$ small enough, and negative for $t$ big enough. Hence, by Bolzano's theorem, there exists some $t=t_\lambda$ such that $\varphi_1(t_\lambda u)=0$. (This
$t_\lambda$ needs not to be unique, but this does not matter for our purposes).

In order to give an upper bound for $t_\lambda$, it is enough to find some $t_1$, such that $\varphi_1(t_1 w_0)<0$. We observe that:
$$
\varphi_1(t w_0) \leq \max\{A t^{p^-} -  \lambda c_3 C t^{r^+};A t^{p^+} -  \lambda c_3 C t^{r^-}\}
$$
so it is enough to choose $t_1$ such that $\max\{A t_1^{p^-} -  \lambda c_3 C t_1^{r^+};A t_1^{p^+} -  \lambda c_3 C t_1^{r^-}\}=0$, i.e.,
$$
t_1 = \left(\frac{A}{c_3 \lambda C}\right)^{1/(r^+-p^-)}\mbox{ or }t_1 = \left(\frac{A}{c_3 \lambda C}\right)^{1/(r^--p^+)}.
$$
Hence, again by Bolzano's theorem, we can choose $t_\lambda \in [0,t_1]$, which
implies that $t_\lambda \to 0$ as $\lambda \to +\infty$.
\end{proof}

For the proof of the Theorem, we need also the following Lemmas.

\begin{lema}\label{lema1}
There exists $C_1, C_2>0$ depending on $p(x)$ and on $c_2$ such that, for every $u\in K_i$, $i=1,2,3$, it holds
$$
\int_\Omega |\nabla u|^{p(x)}\, dx = \left(\lambda \int_\Omega f(x,u)u\, dx
+ \int_\Omega |u|^{q(x)}\, dx\right) \leq C_1 \Phi(u) \leq C_2 \left(\int_\Omega
|\nabla u|^{p(x)}\, dx\right).
$$
\end{lema}

\begin{proof}
The equality is clear since $u\in K_i$.

Now, by (F2), $F(x,u)\ge 0$ so
\begin{align*}
\Phi(u) &= \int_{\Omega} \frac{1}{p(x)} |\nabla u|^{p(x)} - \frac{1}{q(x)} |u|^{q(x)} -
\lambda F(x,u)\, dx \\
&\leq \frac{1}{p^-} \int_\Omega |\nabla u|^{p(x)}\, dx.
\end{align*}

To prove final inequality we proceed as follows, using the norming condition of $K_i$ and hypothesis (F2):
\begin{align*}
\Phi(u)=&\int_{\Omega}\frac{1}{p(x)}|\nabla u|^{p(x)}-\frac{1}{q(x)}|u|^{q(x)}-\lambda
F(x,u)dx\\
& \geq\left(\frac{1}{p^+}-\frac{1}{q^-}\right) \int_{\Omega} |u|^{q(x)}\, dx + \lambda \int_{\Omega} \left(\frac{1}{p^+} f(x,u)u - F(x,u)\right)\, dx\\
&\geq \left(\frac{1}{p^+}-\frac{1}{q^-}\right) \int_{\Omega} |u|^{q(x)}\, dx + \left(\frac{1}{p^+}-\frac{1}{c_2}\right)\lambda\int_\Omega f(x,u)u dx.
\end{align*}

(Recall that $q^->p^+$).This finishes the proof.
\end{proof}

\begin{lema}\label{lema2}
There exists $c>0$ such that
\begin{align*}
\|\nabla u_+\|_{L^{^{p(x)}}(\Omega)} & \geq c\quad\forall\,u\in K_1,\\
\|\nabla u_-\|_{L^{^{p(x)}}(\Omega)} & \geq c\quad\forall\,u\,\in K_2,\\
\|\nabla u_+\|_{L^{^{p(x)}}(\Omega)}\, \mbox{,}\,\|\nabla u_-\|_{L^{^{p(x)}}(\Omega)} &
\geq c\quad\forall\,u\in K_3.
\end{align*}
\end{lema}

\begin{proof}
Suppose that $\|\nabla u_\pm\|_{L^{^{p(x)}}(\Omega)}<1$ By the definition of $K_i$,by (F2) and the Poincar\'e inequality we have that
\begin{align*}
\|\nabla u_\pm\|^{p^+}_{L^{p(x)}(\Omega)}&\leq\rho_{p}(\nabla u_\pm)=\int_\Omega\lambda f(x,u)u_\pm+|u_\pm|^{q(x)} dx\\
&\leq C\rho_{r}(u_\pm)+\rho_{q}(u_\pm)\\
&\leq C\| u_\pm\|^{r^-}_{L^{r(x)}(\Omega)}+\|u_\pm\|^{q^-}_{L^{q(x)}(\Omega)}\\
&\leq c_1\| \nabla u_\pm\|^{r^-}_{L^{p(x)}(\Omega)}+ c_2\| \nabla u_\pm\|^{q^-}_{L^{p(x)}(\Omega)}.
\end{align*}
As $p^+<r^-<q^-$, this finishes the proof.
\end{proof}



The following lemma describes the properties of the manifolds $M_i$.

\begin{lema}\label{lema4}
$M_i$ is a $C^{1}$ sub-manifold of $W_0^{1,p(x)}(\Omega)$ of co-dimension 1 $(i=1,2)$, 2 $(i=3)$ respectively. The sets $K_i$ are complete. Moreover, for every $u\in M_i$ we have the direct decomposition
$$
T_u W_0^{1,p(x)}(\Omega) = T_u M_i \oplus \mbox{span}\{u_+, u_-\},
$$
where $T_u M$ is the tangent space at $u$ of the Banach manifold $M$. Finally, the projection onto the first component in this decomposition is uniformly continuous on bounded sets of $M_i$.
\end{lema}

\begin{proof}
Let us denote
\begin{align*}
&\bar M_1 = \left\{u\in W^{1,p(x)}_0(\Omega)\colon \int_\Omega u_+\, dx > 0
\right\},\\
&\bar M_2 = \left\{u\in W^{1,p(x)}_0(\Omega)\colon \int_\Omega u_-\, dx > 0
\right\},\\
&\bar M_3 = \bar M_1 \cap \bar M_2.
\end{align*}
Observe that $M_i\subset \bar M_i$.

The set $\bar M_i$ is open in $W^{1,p(x)}(\Omega)$,
therefore it is enough to prove that $M_i$ is a $C^1$ sub-manifold of $\bar M_i$. In order to do this, we will construct a $C^{1}$ function $\varphi_i:\bar M_i\to \R^d$ with $d=1\ (i=1,2)$, $d=2\ (i=3)$ respectively and
$M_i$ will be the inverse image of a regular value of $\varphi_i$.

In fact, we define: For $u\in \bar M_1,$
$$
\varphi_1(u)=\int_{\Omega}|\nabla u_+|^{p(x)}-|u_+|^{q(x)}-\lambda f(x,u) u_+ \,
dx.
$$
For $u\in \bar M_2,$
$$
\varphi_2(u)=\int_{\Omega}|\nabla u_-|^{p(x)}-|u_-|^{q(x)}-\lambda f(x,u) u_-\, dx.
$$
For $u\in \bar M_3,$
$$
\varphi_3(u)=(\varphi_1(u),\varphi_2(u)).
$$
Obviously, we have $M_i = \varphi_i^{-1}(0)$. From standard arguments (see \cite{DJM}, or the appendix of \cite{R}), $\varphi_i$ is of class $C^1$. Therefore, we only need to show that $0$ is a regular value for $\varphi_i$. To this end we compute, for $u\in M_1$,
\begin{align*}
\langle\nabla\varphi_1(u),u_+\rangle&\leq p^+ \rho_{p}(\nabla u_+)-q^-\rho_q(u_+)-\lambda\int_\Omega f(x,u)u_+- f_u(x,u)u_+^2\,dx\\
&\leq q^-\left(\rho_p(\nabla u_+)-\rho_q(u_+)\right)-\lambda\int_\Omega f(x,u)u_+- f_u(x,u)u_+^2\,dx\\
&\leq (q^-\lambda-\lambda)\int_\Omega f(x,u)u_+\,dx-\int_\Omega
f_u(x,u)u_+^2\,dx.
\end{align*}
By (F2) the last term is bounded by
\begin{align*}
(q^-\lambda-\lambda-\frac{\lambda}{c_1})\int_\Omega f(x,u)u_+\,dx
&=\left(q^--1-\frac{1}{c_1}\right)\left(\rho_p(\nabla u_+)-\rho_q(u_+)\right)\\
&\leq\left(q^--1-\frac{1}{c_1}\right)\rho_p(\nabla u_+).
\end{align*}
Recall that $c_1<1/(q^--1)$. Now, the last term is strictly negative by Lemma \ref{lema2}. Therefore, $M_1$ is a $C^{1}$ sub-manifold of $W^{1,p(x)}(\Omega)$.
The exact same argument applies to $M_2$. Since trivially
$$
\langle \nabla \varphi_1(u), u_-\rangle = \langle \nabla \varphi_2(u),
u_+\rangle =0
$$
for $u\in M_3$, the same conclusion holds for $M_3$.

To see that $K_i$ is complete, let $u_k$ be a Cauchy sequence in $K_i$, then $u_k\to u$ in $W^{1,p(x)}(\Omega)$. Moreover, $(u_k)_{\pm}\to u_{\pm}$ in $W^{1,p(x)}(\Omega)$. Now it is easy to see, by Lemma \ref{lema2} and by continuity that $u\in K_i$.

Finally, by the first part of the proof we have the decomposition
$$
T_u W^{1,p(x)}(\Omega) = T_u M_i \oplus \mbox{span}\{u_+\}
$$
Where $M_1=\{u:\varphi_1(u)=0\}$ and $T_u M_1=\{v:\langle\nabla \varphi_1(u),v\rangle=0\}$. Now let $v\in T_u W_0^{1,p(x)}(\Omega)$ be a unit tangential vector, then $v = v_1 + v_2$ where $v_2=\alpha u_+$ and $v_1=v-v_2$. Let us take $\alpha$ as
$$
\alpha = \frac{\langle \nabla\varphi_1(u),v\rangle}{\langle
\nabla\varphi_1(u),u_+\rangle}.
$$
With this choice, we have that $v_1\in T_u M_1$. Now
$$
\langle\varphi_1(u),v_1\rangle=0.
$$
The very same argument to show that $T_uW^{1,p(x)}(\Omega) = T_u M_2\oplus\langle u_-\rangle$ and $T_uW^{1,p(x)}(\Omega) = T_u M_3\oplus \langle u_+,u_-\rangle$.

From these formulas and from the estimates given in the first part of the proof, the uniform continuity of the projections onto $T_u M_i$ follows.
\end{proof}

Now, we need to check the Palais-Smale condition for the functional $\Phi$ restricted to the manifold $M_i$.
We begin by proving the Palais-Smale condition for the functional $\Phi$ unrestricted, below certain level of energy.

\begin{lema}\label{acotada}
Assume that $r\le q$. Let $\{u_j\}_{j\in\N}\subset W_0^{1,p(x)}(\Omega)$ a Palais-Smale sequence then $\{u_j\}_{j\in\N}$ is bounded in $W_0^{1,p(x)}(\Omega)$.
\end{lema}
\begin{proof}
By definition
$$
\Phi(u_j)\to c \qquad \mbox{and}\qquad\Phi'(u_j)\to 0.
$$
Now, we have
$$
c+1 \geq \Phi(u_j) = \Phi(u_j) - \frac{1}{c_2} \langle \Phi'(u_j), u_j \rangle + \frac{1}{c_2} \langle \Phi'(u_j), u_j \rangle,
$$
where
$$
\langle \Phi'(u_j), u_j \rangle = \int_\Omega |\nabla u_j|^{p(x)} - |u_j|^{q(x)} -\lambda f(x,u_j)u_j\, dx.
$$
Then, if $c_2<q^-$ we conclude
$$
c+1 \geq \left(\frac{1}{p+} - \frac{1}{c_2}\right) \int_\Omega |\nabla u_j|^{p(x)}\, dx - \frac{1}{c_2} |\langle \mathcal{F}'(u_j), u_j \rangle|.
$$
We can assume that $\|u_j\|_{W_0^{1,p(x)}(\Omega)}\geq 1$. As $\|\mathcal{F}'(u_j)\|$ is bounded we have that
$$
c+1 \geq \left(\frac{1}{p+} - \frac{1}{c_2}\right) \|u_j\|^{p^-}_{W_0^{1,p(x)}(\Omega)} - \frac{C}{c_2}\|u_j\|_{W_0^{1,p(x)}(\Omega)}.
$$
We deduce that $u_j$ is bounded.

This finishes the proof.
\end{proof}

From the fact that $\{u_j\}_{j\in\N}$ is a Palais-Smale sequence it follows, by Lemma \ref{acotada}, that $\{u_j\}_{j\in\N}$ is bounded in $W_0^{1,p(x)}(\Omega)$. Hence, by The Concentration-Compactness method for variable exponent (See\cite{FBS}), we have
\begin{align}
&|u_j|^{q(x)}\rightharpoonup \nu=|u|^{q(x)} + \sum_{i\in I} \nu_i\delta_{x_i} \quad \nu_i>0,\\
&|\nabla u_j|^{p(x)}\rightharpoonup \mu \geq |\nabla u|^{p(x)} + \sum_{i\in I} \mu_i \delta_{x_i}\quad \mu_i>0,\\
&S \nu_i^{1/p^*(x_i)} \leq \mu_i^{1/p(x_i)}.
\end{align}

Note that if $I=\emptyset$ then $u_j\to u$ strongly in $L^{q(x)}(\Omega)$. We know that $\{x_i\}_{i\in I}\subset \mathcal{A}:=\{x:q(x)=p^*(x)\}$.We define $q^-_{\mathcal A}:=\inf_{\mathcal A} q(x).$

Let us show that if $c < \left(\frac{1}{p^+} - \frac{1}{q^-_\mathcal{A}}\right)S^n$ and $\{u_j\}_{j\in\N}$ is a Palais-Smale sequence, with energy level $c$, then $I=\emptyset$.

\medskip

In fact, suppose that $I \not= \emptyset$. Then let $\phi\in C_0^\infty(\R^n)$ with support in the unit ball of $\R^n$. Consider the rescaled functions $\phi_{i,\varepsilon}(x) = \phi(\frac{x-x_i}{\varepsilon})$.

As $\Phi'(u_j)\to 0$ in $(W_0^{1,p(x)}(\Omega))'$, we obtain that
$$
\lim_{j\to\infty} \langle \Phi'(u_j), \phi_{i,\varepsilon}u_j \rangle = 0.
$$
On the other hand,
$$
\langle \Phi'(u_j), \phi_{i,\varepsilon} u_j \rangle = \int_\Omega |\nabla u_j|^{p(x)-2}\nabla u_j\nabla(\phi_{i,\varepsilon}u_j) - \lambda f(x,u_j)u_j\phi_{i,\varepsilon} - |u_j|^{q(x)}\phi_{i,\varepsilon}\, dx
$$
Then, passing to the limit as $j\to\infty$, we get
\begin{align*}
0 = &\lim_{j\to\infty} \left(\int_\Omega |\nabla u_j|^{p(x)-2} \nabla u_j \nabla(\phi_{i,\varepsilon}) u_j\, dx\right)\\
& + \int_\Omega \phi_{i,\varepsilon}\, d\mu - \int_\Omega \phi_{i,\varepsilon}\, d\nu
- \int_\Omega\lambda f(x,u)u\phi_{i,\varepsilon}\, dx.
\end{align*}
By H\"older inequality, it is easy to check that
$$
\lim_{j\to\infty} \int_\Omega |\nabla u_j|^{p(x)-2} \nabla u_j \nabla(\phi_{i,\varepsilon}) u_j\, dx = 0.
$$
On the other hand,
$$
\lim_{\varepsilon\to 0} \int_\Omega \phi_{i,\varepsilon}\, d\mu = \mu_i\phi(0),\qquad
\lim_{\varepsilon\to 0} \int_\Omega \phi_{i,\varepsilon}\, d\nu = \nu_i\phi(0).
$$
and
$$
\lim_{\varepsilon\to 0} \int_\Omega\lambda f(x,u)u\phi_{i,\varepsilon}\, dx = 0.
$$
So, we conclude that $(\mu_i-\nu_i)\phi(0)=0$, i.e, $\mu_i=\nu_i$. Then,
$$
S \nu_i^{1/p^*(x_i)} \leq \nu_i^{1/p(x_i)},
$$
so it is clear that $\nu_i = 0$ or $S^n\leq\nu_i$.

On the other hand, as $c_2>p^+$,
\begin{align*}
c =& \lim_{j\to\infty} \Phi(u_j) = \lim_{j\to\infty} \Phi(u_j) - \frac{1}{p+} \langle \Phi'(u_j), u_j \rangle\\
=& \lim_{j\to\infty} \int_\Omega \left(\frac{1}{p(x)} - \frac{1}{p+}\right) |\nabla u_j|^{p(x)}\, dx + \int_\Omega \left(\frac{1}{p+} - \frac{1}{q(x)}\right) |u_j|^{q(x)}\, dx\\
& - \lambda \int_\Omega F(x,u_j)\, dx+\frac{\lambda}{p^+}\int_\Omega f(x,u_j)u_j\,dx\\
\geq& \lim_{j\to\infty} \int_\Omega \left(\frac{1}{p+}-\frac{1}{q(x)}\right) |u_j|^{q(x)}\, dx\\
\geq& \lim_{j\to\infty} \int_{\mathcal{A}_\delta} \left(\frac{1}{p+}-\frac{1}{q(x)}\right) |u_j|^{q(x)}\, dx\\
\geq&\lim_{j\to\infty} \int_{\mathcal{A}_\delta} \left(\frac{1}{p+}-\frac{1}{q^-_{\mathcal{A}_\delta}}\right) |u_j|^{q(x)}\, dx
\end{align*}
But
\begin{align*}
\lim_{j\to\infty} \int_{\mathcal{A}_\delta} \left(\frac{1}{p+} - \frac{1}{q^-_{\mathcal{A}_\delta}}\right) |u_j|^{q(x)}\, dx &= \left(\frac{1}{p+}-\frac{1}{q^-_{\mathcal{A}_\delta}}\right) \left(\int_{\mathcal{A}_\delta}|u|^{q(x)}\, dx + \sum_{j\in I} \nu_j\right)\\
&\geq \left(\frac{1}{p+} - \frac{1}{q^-_{\mathcal{A}_\delta}}\right) \nu_i\\
&\geq \left(\frac{1}{p+}-\frac{1}{q^-_{\mathcal{A}_\delta}}\right) S^n.
\end{align*}
As $\delta>0$ is arbitrary, and $q$ is continuous, we get
$$
c\ge \left(\frac{1}{p+}-\frac{1}{q^-_{\mathcal{A}}}\right) S^n.
$$

Therefore, if
$$
c < \left(\frac{1}{p+} - \frac{1}{q^-_{\mathcal{A}}}\right)S^n,
$$
the index set $I$ is empty.

Now we are ready to prove the Palais-Smale condition below level $c$.
\begin{lema}\label{PS sin restringir}
Let $\{u_j\}_{j\in\N}\subset W_0^{1,p(x)}(\Omega)$ be a  Palais-Smale sequence, with energy level $c$. If $c < \left(\frac{1}{p+} - \frac{1}{q^-_{\mathcal{A}}}\right) S^n$, then there exist $u\in W_0^{1,p(x)}(\Omega)$ and $\{u_{j_k}\}_{k\in\N}\subset \{u_j\}_{j\in\N}$ a subsequence such that $u_{j_{k}}\to u$ strongly in $W_0^{1,p(x)}(\Omega)$.
\end{lema}

\begin{proof}
We have that $\{u_j\}_{j\in\N}$ is bounded. Then, for a subsequence that we still denote $\{u_j\}_{j\in\N}$, $u_j\to u$ strongly in $L^{q(x)}(\Omega)$. We define $\Phi'(u_j):=\phi_j$. By the Palais-Smale condition, with energy level c, we have $\phi_j\to 0$ in $(W_0^{1,p(x)}(\Omega))'$.

By definition $\langle \Phi'(u_j), z \rangle = \langle \phi_j, z \rangle$ for all $z\in W_0^{1,p(x)}(\Omega)$, i.e,
$$
\int_\Omega |\nabla u_j|^{p(x)-2}\nabla u_j\nabla z\, dx - \int_\Omega |u_j|^{q(x)-2} u_j z\, dx - \int_\Omega\lambda f(x,u_j) z\, dx = \langle \phi_j, z \rangle.
$$
Then, $u_j$ is a weak solution of the following equation.
\begin{equation}
\begin{cases}
-\Delta_{p(x)}u_j=|u_j|^{q(x)-2}u_j+\lambda f(x,u_j)+\phi_j=: f_j &\mbox{in }\Omega,\\
u_j = 0 &\mbox{on }\partial\Omega.
\end{cases}
\end{equation}
We define $T\colon (W_0^{1,p(x)}(\Omega))' \to W_0^{1,p(x)}(\Omega)$, $T(f):=u$ where $u$ is the weak solution of the following equation.
\begin{equation}
\begin{cases}
-\Delta_{p(x)}u=f & \mbox{in } \Omega,\\
u = 0 & \mbox{on } \partial\Omega.
\end{cases}
\end{equation}
Then $T$ is a continuous invertible operator.

It is sufficient to show that $f_j$ converges in $(W^{1,p(x)}_0(\Omega))'$. We only need to prove that $|u_j|^{q(x)-2}u_j \to |u|^{q(x)-2}u$ strongly in $(W_0^{1,p(x)}(\Omega))'$.

In fact,
\begin{align*}
\langle|u_j|^{q(x)-2}u_j-|u|^{q(x)-2}u,\psi\rangle&=\int_\Omega(|u_j|^{q(x)-2}u_j-|u|^{q(x)-2}u)\psi\,dx\\
&\leq\|\psi\|_{L^{q(x)}(\Omega)}\|(|u_j|^{q(x)-2}u_j-|u|^{q(x)-2}u)\|_{L^{q'(x)}(\Omega)}.
\end{align*}
Therefore,
\begin{align*}
\|(|u_j|^{q(x)-2}u_j - |u|^{q(x)-2}u)\|_{(W_0^{1,p(x)}(\Omega))'} &= \sup_{\genfrac{}{}{0cm}{}{\psi\in W^{1,p(x)}_0(\Omega)} {\|\psi\|_{W^{1,p(x)}_0(\Omega)}=1}} \int_\Omega (|u_j|^{q(x)-2}u_j-|u|^{q(x)-2}u)\psi\, dx\\
&\leq \|(|u_j|^{q(x)-2}u_j - |u|^{q(x)-2}u)\|_{L^{q'(x)}(\Omega)}
\end{align*}
and now, by the Dominated Convergence Theorem this last term goes to zero as $j\to\infty$.

The proof is finished.
\end{proof}



Now, we can prove the Palais-Smale condition for the restricted functional.

\begin{lema}\label{lema5}
The functional $\Phi|_{K_i}$ satisfies the Palais-Smale condition for energy level $c$ for every
$c < \left(\frac{1}{p+}-\frac{1}{q^-_{\mathcal{A}}}\right)S^n$.
\end{lema}

\begin{proof}
Let $\{u_k\}\subset K_i$ be a Palais-Smale sequence, that is $\Phi(u_k)$ is uniformly bounded and $\nabla \Phi|_{K_i}(u_k)\to 0$ strongly. We need to show that there exists a subsequence $u_{k_j}$ that converges strongly in $K_i$.

Let $v_j\in T_{u_j}W_0^{1,p(x)}(\Omega)$ be a unit tangential vector such that
$$
\langle \nabla \Phi(u_j), v_j\rangle = \|\nabla
\Phi(u_j)\|_{(W^{1,p(x)}_0(\Omega))'}.
$$
Now, by Lemma \ref{lema4}, $v_j = w_j + z_j$ with $w_j\in
T_{u_j}M_i$ and $z_j\in \mbox{span}\{(u_j)_+, (u_j)_-\}$.

Since $\Phi(u_j)$ is uniformly bounded, by Lemma \ref{lema1}, $u_j$ is uniformly bounded in $W_0^{1,p(x)}(\Omega)$ and hence $w_j$ is uniformly bounded
in $W_0^{1,p(x)}(\Omega)$. Therefore
$$
\|\nabla\Phi(u_j)\|_{(W^{1,p(x)}_0(\Omega))'} = \langle \nabla \Phi(u_j), v_j\rangle
= \langle \nabla \Phi|_{K_i}(u_j), v_j\rangle\to 0.
$$

As $w_j$ is uniformly bounded and $\nabla \Phi|_{K_i}(u_k)\to 0$ strongly, the inequality converges strongly to 0. Now the result follows by Lema \ref{PS sin restringir}.
\end{proof}

We now immediately obtain

\begin{lema}\label{lema6}
Let $u\in K_i$ be a critical point of the restricted functional $\Phi|_{K_i}$. Then $u$ is also a critical point of the unrestricted functional $\Phi$ and hence a weak solution to \eqref{1.1}.
\end{lema}

With all this preparatives, the proof of the Theorem follows easily.
\begin{proof}[\bf Proof of Theorem 1]
To prove the Theorem,  we need to check that the functional $\Phi\mid_{K_i}$ verifies the hypotheses of the Ekeland's Variational Principle \cite{Ekeland}.

The fact that $\Phi$ is bounded below over $K_i$ is a direct consequence of the construction of the manifold $K_i$.

Then, by Ekeland's Variational Principle, there existe $v_k\in K_i$ such that
$$
\Phi(v_k)\to c_i:=\inf_{K_i}\Phi\qquad \mbox{ and }\qquad
(\Phi\mid_{K_i})'(v_k)\to 0.
$$
We have to check that if we choose $\lambda$ large, we have that $c_i< \left(\frac{1}{p+} - \frac{1}{q^-_{\mathcal{A}}}\right)S^n,$. This follows easily from Lemma \ref{tlambda}. For instance, for $c_1$, we have that choosing $w_0\geq 0$, if $t_\lambda<1$
$$
c_1 \leq \Phi(t_\lambda w_0) \leq \frac{1}{p^-} t_{\lambda}^{p^+} \int_\Omega |\nabla w_0|^{p(x)}
\, dx
$$
Hence $c_1 \to 0$ as $\lambda \to \infty$. Moreover, it follows from the estimate of $t_\lambda$ in Lemma \ref{tlambda}, that $c_i<\left(\frac{1}{p+} - \frac{1}{q^-_{\mathcal{A}}}\right)S^n$ for $\lambda>\lambda^*(p,q,n,c_3)$. The other cases are similar.

From Lemma \ref{lema5}, it follows that $v_k$ has a convergent subsequence, that we still call $v_k$. Therefore $\Phi$ has a critical point in $K_i$, $i=1,2,3$
and, by construction, one of them is positive, other is negative and the last one changes sign.
\end{proof}

\textbf{Acknowledgements} I want to thank Juli\'an Fern\'andez Bonder for valuable help.

\end{document}